\documentclass[11pt]{amsart}
\usepackage{fullpage}
\usepackage{amsfonts}
\usepackage{latexsym}
\usepackage{amscd}
\usepackage{amsmath}
\usepackage{amssymb}
\usepackage{graphicx}
\usepackage[mathscr]{eucal}
\usepackage{xypic}

\newcommand{\finp }{\hspace{\fill}\rule{.2cm}{.2cm}\medskip}

\newcommand{\demost}{\noindent{\itshape\bfseries Proof.} }

\newtheorem{teorema}{Theorem}[section]
\newtheorem{lema}[teorema]{Lemma}
\newtheorem{corolario}[teorema]{Corollary}
\newtheorem{proposicion}[teorema]{Proposition}

\theoremstyle{definition}

\newtheorem{example}[teorema]{Example}

\theoremstyle{remark}

\numberwithin{equation}{section}
\newcommand{\N}{\mathbb{N}}

\begin{document}

\title{Iterates of dynamical systems on compact  metrizable countable spaces}

\author{S. Garc\'{\i}a-Ferreira}
\address{Centro de Ciencias Matem\'aticas, Universidad Nacional
Aut\'onoma de M\'exico, Apartado Postal 61-3, Santa Mar\'{\i}a,
58089, Morelia, Michoac\'an, M\'exico}
%\curraddr{Department of Mathematics and Statistics,
%Case Western Reserve University, Cleveland, Ohio 43403}
\email{sgarcia@matmor.unam.mx}

\author{Y. Rodriguez-L\'opez}
\address{Secci\'on de Matem\'aticas, Universidad Nacional Experimental Polit\'ecnica ``Antonio Jose de Sucre", Barquisimeto, Venezuela}
\email{yrodriguez@unexpo.edu.ve}

\author{C. Uzc\'ategui}
\address{Departamento de Matem\'{a}ticas,
Facultad de Ciencias, Universidad de los Andes, M\'erida 5101,
Venezuela} \email{uzca@ula.ve}

\thanks{The research of the first listed author was supported by  CONACYT grant no. 176202 and PAPIIT
grant no. IN-101911. The research of the second listed author was
supported by Universidad Nacional Experimental Polit\'ecnica
``Antonio Jose de Sucre", Venezuela. Hospitality and financial
support received from the \emph{Universidad de Los Andes,
Venezuela} and the  {\em Centro de Ciencias Matem\'aticas de la
Universidad Aut\'onoma de M\'exico (Morelia)} where this research
was performed are gratefully acknowledged}

%    General info

\date{}

\dedicatory{}


\begin{abstract}
Given a dynamical system $(X,f)$, we let $E(X,f)$ denote its Ellis
semigroup and $E(X,f)^* = E(X,f) \setminus \{ f^n : n \in
\mathbb{N}\}$. We analyze the Ellis semigroup of a dynamical
system having a compact metric countable space as a phase space.
We show that if  $(X,f)$ is a dynamical system such that $X$ is a
compact metric countable space and  every accumulation point of
$X$ is periodic, then either all function of $E(X,f)^*$ are
continuous or all functions of $E(X,f)^*$ are discontinuous. We
describe an example of a  dynamical system $(X,f)$ where $X$ is a
compact metric countable space,  the orbit of each accumulation
point is finite and $E(X,f)^*$ contains both continuous and
discontinuous functions.
\end{abstract}

\maketitle

%%%%%%%%%%%%%%%%%%%%%%%%%%%%%%%%%%%%%%%%%%%%%%%%%%%%%%%%%%%%
\section{Introduction}
%%%%%%%%%%%%%%%%%%%%%%%%%%%%%%%

 We start the paper by fixing some standard notions and terminology. Let $(X,f)$ be a dynamical system.
 The {\it orbit} of $x$, denoted by $\mathcal O_f(x)$, is
the set $\{f^n(x):n\in\mathbb N\}$, where $f^n$ is $f$ composed
with itself $n$ times. A point $x\in X$ is called a {\it periodic
point} of $f$ if there exists $n\geq 1$ such that $f^n(x)=x$, and
$x$ is called {\it eventually periodic} if its orbit is finite.
The $\omega-${\it limit set} of $x\in X$, denoted by
$\omega_f(x)$, is the set of points $y\in X$ for which there
exists an increasing sequence $(n_k)_{k\in\mathbb N}$ such that
$f^{n_k}(x)\rightarrow y$. For each $y\in\mathcal O_f(x)$,
$\omega_f(y)=\omega_f(x)$. If $\mathcal{O}_f(y)$ contains a
periodic point $x$, then $\omega_f(y)=\mathcal{O}_f(x)$. We denote
by $\mathcal{N}(x)$ the collection of all the neighborhoods of
$x$, for each $x\in X$. The set of all accumulation points of $X$,
the derivative of $X$, is denoted by $X'$. We remark that the
countable ordinal space $\omega^2+1$ is homeomorphic to the
compact metric subspace $Y=\{ 1 - \frac{1}{n}:n\in\mathbb N
\setminus\{0\}\}\cup\{1\}\cup (\bigcup_{n\in\mathbb N}A_n)$ of
$\mathbb{R}$, where $A_n$ is an increasing sequence contained in
$(1-\frac{1}{n-1},1-\frac{1}{n})$ such that $A_n\longrightarrow
1-\frac{1}{n}$, for each $n \in \mathbb{N}$ bigger than $1$. The
Stone-\v{C}ech compactification $\beta(\mathbb N)$ of  $\mathbb N$
with the discrete topology will be identified with the set of
ultrafilters over $\mathbb N$. Its remainder ${\mathbb N}^*=
\beta({\mathbb N})\setminus \mathbb{N}$ is  the set of all free
ultrafilters on $\mathbb N$, where, as usual, each natural number
$n$ is identified with  the fixed ultrafilter consisting of all
subsets of $\mathbb N$ containing $n$. For $A\subseteq \mathbb N$,
$A^*$ denotes the collection of all $p\in \mathbb{N}^*$ such that
$A\in p$.

\medskip

In our dynamical systems $(X,f)$ the space $X$ will be  compact
metric  and $f:X\rightarrow X$  will be a continuous map. A very
useful object to study the topological behavior of the dynamical
system $(X,f)$ is the so-called  {\em Ellis semigroup} or {\em
enveloping semigroup}, introduced by Ellis \cite{ell}, which is
defined as  the pointwise closure of $\{f^n:\; n\in \mathbb N\}$
in the compact space $X^X$ with composition of functions as its
algebraic operation.  The Ellis semigroup, denoted $E(X,f)$, is
equipped with the topology inhered from the product space $X^X.$
Enveloping semigroups have played a very crucial role in
topological dynamics and they are an active area of research (see,
for instance, the survey article \cite{glasner}).

\medskip

The motivation of our work is the fact that for some spaces either
all functions of $E^*(X,f)$ are continuous or all  are
discontinuous. Namely, this holds when $X$ is a convergent
sequence with its limit point \cite{GarciaSanchis} (see also
\cite{Garcia}) and for $X=[0,1]$ as it was recently shown by P.
Szuca \cite{Szuca}. In this direction, we will show that it also
happens for any dynamical system $(X,f)$ where $X$ is a compact
metrizable countable space such that every accumulation point of
$X$ is periodic.  We also present an example of a dynamical system
$(X,f)$ where $X$ is the ordinal space $\omega^2+1$ such that
$E(X,f)\setminus \{f^n:\; n\in\mathbb N\}$ contains continuous and
also discontinuous functions and each accumulation point of $X$ is
eventually periodic. This answers a question posed in
\cite{GarciaSanchis}.

\medskip

Now we recall a convenient description of $E(X,f)$  in terms of
the notion of $p-$limits where $p$  is an ultrafilter on the
natural number $\N$. Given $p\in \mathbb{N}^*$ and a sequence
$(x_n)_{n\in\mathbb{N}}$ in a space $X$, we say that a point $x\in
X$ is the $p-${\it limit point} of the sequence, in symbols $x
=p-\lim_{n\rightarrow \infty}x_n$, if for every neighborhood $V$
of $x$, $\{n\in\mathbb{N}: f^n(x)\in V\} \in p$. Observe that a
point $x\in X$ is an accumulation point of a countable set
$\{x_n:\,n\in\mathbb{N}\}$ of $X$ iff there is $p\in \mathbb{N}^*$
such that $x = p-\lim_{n\rightarrow \infty}x_n$. It is not hard to
prove that each sequence of a compact space always has a p-limit
point for every $p\in \mathbb{N}^*$. The notion of a $p-$limit
point has been used in topology and analysis (see for instance
\cite{Be} and \cite[p. 179]{fu}).

\medskip

A. Blass \cite{Bla} and N. Hindman \cite{hi} formally established
the connection between  ``the iteration in to\-po\-lo\-gi\-cal
dynamics'' and ``the convergence with respect to an ultrafilter''
by considering a more general iteration of the function $f$ as
follows: Let $X$ be
 a compact  space and $f : X
\rightarrow  X$ a continuous function. For $p\in\mathbb{N}^*$, the
$p-$iterate of $f$ is the function $f^p: X\rightarrow X$ defined
by \[ f^p(x) = p-\lim_{n\rightarrow \infty} f^n(x),
\]
for all $x\in X$. The description of the Ellis semigroup in terms
of the $p-$iterates is then the following:
\[
\begin{array}{rcl}
E(X,f) & = &\{f^p: p\in \beta\mathbb N \}\\
\\
f^p\circ f^q & =& f^{q+p}\;\; \mbox{for each $p,q\in\beta\mathbb
N$ (see \cite{Bla}, \cite{hi})}.
\end{array}
\]

\medskip

The paper is organized as follows. The second section is devoted
to prove some basic results that will be  used in the rest of the
paper. In the third section, we show our main results about
$E(X,f)$ when $X$ is a compact metric countable space and each
element of $X'$ is a periodic point of $f$. In the forth section,
we construct a dynamical system $(X,f)$ in which all accumulation
points are eventually periodic and $E(X,f) \setminus \{ f^n : n
\in \mathbb{N}\}$ contains continuous and also discontinuous
functions. We also state some open questions.

%%%%%%%%%%%%%%%%%%%%%%%%%%%%%%%%%%%%%%%%%%%%%%%%%%%%%%%%%%%%%%
\section{Basic properties}
%%%%%%%%%%%%%%%%%%%%%%%%%%%%%%%%%%%%%%%%%%%%%%%%%%%%%%%%%%%%%%

In this section, we present some basic lemmas that will be used in
the sequel.

\begin{lema}
\label{iteradas de potencias} Let $(X,f)$ be a dynamical system
and $p\in \mathbb{N}^*$. If $g=f^n$ for some positive $n \in \N$,
then
$$
g^p=f^p\circ f^n.
$$
\end{lema}

\proof By definition, we have that
$$
 g^p(x) = p-\lim_{k \rightarrow \infty}g^k(x) = p-\lim_{k \rightarrow \infty}f^n(f^k(x))
$$
$$
= f^n\big(p-\lim_{k \rightarrow \infty}f^k(x)\big) = f^n\circ
f^p(x) = f^p\circ f^n(x),
$$
for every $x \in X$. \finp

Now, we calculate the $p$-iteration at certain points of a
dynamical system.

\begin{lema}\label{piteradaperiodico2} Let $(X,f)$ be a dynamical system and let $x \in X$ be a periodic point. If
$x$ has period $n$ and $p \in \big(n \mathbb{N} + l\big)^*$ for
some $l < n$, then $f^p(x) = f^l(x)$.
\end{lema}

\demost Let $V$ be an open neighborhood of $f^p(x)$. By
definition, we have that $\{ k \in \mathbb{N} : f^k(x) \in V \}
\in p$. Thus, $A = (n\mathbb{N} + l) \cap \{ k \in \mathbb{N} :
f^k(x) \in V \} \in p$. For each $k \in A$ choose $m_k \in
\mathbb{N}$ so that $k = nm_k + l$. Then, $f^k(x) = f^{nm_k +
l}(x) = f^l(f^{nm_k}(x)) = f^l(x)$ for each $k \in A$. Hence,
$f^p(x) = f^l(x)$. \finp

When the point $x$ is eventually periodic, we have the following.

\begin{proposicion}
Let $(X,f)$ be a dynamical system and let $x \in X$  with finite
orbit. If $m \in \N$ is the smallest positive integer such that
$f^m(x)$ is a periodic point with period $n$, then for every  $p
\in \mathbb{N}^*$ there is $l < n$
 such that  $f^p(x) = f^l(f^m(x))$.
\end{proposicion}

\demost It is evident that for every positive integer $k \geq m$
there is $0 \leq l < n$ such that $f^k(x) = f^l(f^m(x))$. Hence,
if $p \in  \mathbb{N}^*$, then there is $l < n$
 such that
$$
 f^p(x) = p-\lim_{k \rightarrow \infty}f^k(x) = f^l(f^m(x)).
$$
\finp

Let $(X,f)$ be a dynamical system and assume that $x \in X$ has
infinite orbit. Then, by using the $p$-iterates, we have that $y
\in \omega_f(x)$  iff there is $p \in \N^*$ such that $f^p(x) =
y$. For the case when $\omega_f(x)$ is finite, we have the
following well-known result (see for instance \cite{BC}).

\begin{lema}
\label{aisladoperiodico} Let $(X,f)$ be a dynamical system.  If
$\omega_f(x)$ is finite, then every point of $\omega_f(x)$ is
periodic. In particular, if $\omega_f(x)$ has a point isolated in
$\overline{\mathcal{O}_f(x)}$, then every point of $\omega_f(x)$
is periodic.
\end{lema}

\demost Fix $y \in \omega_f(x)$. Then, it is clear that
$A=\{n\in\mathbb{N}: f^n(x) = y \}$ is infinite. Hence, if $m, n
\in A$ and $m < n$, then  $y = f^m(x) = f^{n}(x) = f^{n -
m}(f^m(x)) = f^{n - m}(y)$. Therefore, $y$ is periodic. If
$\omega_f(x)$ has a point isolated in
$\overline{\mathcal{O}_f(x)}$, it is evident that $\omega_f(x)$ is
finite. \finp

\begin{corolario}\label{recurrente a periodico}
Let $(X,f)$ be a dynamical system. If $x\in X $ is a recurrent
point  and  there is a point in  $\omega_f(x)$ isolated in
$\overline{\mathcal{O}_f(x)}$, then $x$ is periodic.
\end{corolario}

In the following lemma, we express the orbit $\mathcal O_f(a)$ in
terms of $\mathcal O_g(a)$, where $g$ is an iteration  of the
function  $f$.

\begin{lema}\label{descomposici�n de �rbita}
Let $(X,f)$ be a dynamical system. If $g=f^n$ for some positive
$n\in\mathbb N$, then
$$
\mathcal{O}_f(x)=\mathcal{O}_g(x)\cup
f[\mathcal{O}_g(x)]\cup\ldots\cup f^{n-1}[\mathcal{O}_g(x)],
 $$
for every $x \in X$.
\end{lema}

\demost It is evident that $\mathcal{O}_f(x) \cup
f[\mathcal{O}_g(x)]\cup\ldots\cup
f^{n-1}[\mathcal{O}_g(x)]\subseteq \mathcal O_f(x)$.  Let
$m\in\mathbb N$ and choose $t\in\mathbb N$ and  $0 \leq l< n$ so
that $m=tn+l$.  Then, we have that
$$
f^m(x) = f^{tn +l}(x) = f^{l}(f^{nt}(x)) = f^l(g^t(x)) \in
f^l(\mathcal O_g(x)).
$$
Thus, $\mathcal O_f(x) \subseteq \mathcal{O}_f(x) \cup
f[\mathcal{O}_g(x)]\cup\ldots\cup f^{n-1}[\mathcal{O}_g(x)]$.
Therefore,
$$
\mathcal{O}_f(x)=\mathcal{O}_g(x)\cup
f[\mathcal{O}_g(x)]\cup\ldots\cup f^{n-1}[\mathcal{O}_g(x)].
$$
\finp

Next, we shall analyze when the $\omega$-limit set is equal to the
orbit of a periodic point.

\begin{lema}\label{orbitaomega}
Let $(X,f)$ be a dynamical system and let  $x\in X$ be with
infinite orbit.  If there is  $l \in \N$ such that
$f^{ln}(x)\xrightarrow[n \to \infty]{} y$, then $\omega_f(x) =
\mathcal O_f(y)$ and $y$ has period $l$. Conversely, if
$\omega_f(x) = \mathcal O_f(y)$ and $y$ has period $l$, then
$f^{ln}(x)\xrightarrow[n \to \infty]{} f^i(y)$ for some $i<l$.
\end{lema}

\proof  Suppose $f^{ln}(x)\xrightarrow[n \to \infty]{} y$. Let $z
\in\omega_f(x)$.  Then, there is an increasing sequence $(n_k)_{k
\in \N}$ such that $f^{n_k}(x)\xrightarrow[k \to \infty]{} z$.
Choose  $i<l$ so that $\{n_k:k\in\mathbb N\}\cap (l\mathbb{N}+i)$
is infinite. So, $f^{n_k}(x)=f^i(f^{lt_k}(x))$, for some $t_k \in
\N$, for infinitely many $k$'s. Since $f^{lt_k}(x)\xrightarrow[k
\to \infty]{} y$, we must have that $f^{n_k}(x)\xrightarrow[k \to
\infty]{} f^i(y)$ and hence $f^i(y) = z$.  Therefore, $z \in
\mathcal O_f(y)$. This shows that $\omega_f(x) \subseteq \mathcal
O_f(y)$. Since $y \in \omega_f(x)$, we must have that $\mathcal
O_f(y) \subseteq \omega_f(x)$. Clearly,
$f^{l(n+1)}(x)\xrightarrow[n \to \infty]{} f^l(y)$, thus
$f^l(y)=y$.

Conversely, suppose $\omega_f(x) = \mathcal O_f(y)$ and $y$ has
period $l$. By compactness, there are $i<l$ and $(n_k)_k$
increasing such that $f^{n_kl}(x)\rightarrow f^i(y)$. Let $V$ be
an open set such that $V\cap {\mathcal O}_f(y)=\{f^i(y)\}$. Then,
$A=\{n:f^{nl}(x)\in V\}$ is infinite. Since $\lim_{n\in A}
f^{(n+1)l}(x)= f^i(y)$, then $(A+1)\setminus A$ is finite,
therefore $A$ is a final segment of $\mathbb{N}$ and thus
$f^{ln}(x)\xrightarrow[n \to \infty]{} f^i(y)$. \finp

%%%%%%%%%%%%%%%%%%%%%%%%%%%%%%%%%%%%%%%%%%%%%%%
\section{Main results}
%%%%%%%%%%%%%%%%%%%%%%%%%%%%%%%%%%%%%%%%%%%%%%%

Since our spaces are scattered, the Cantor-Bendixson rank will be
very useful to carry out some inductive process  in several
proofs:

\smallskip

For a successor ordinal $\alpha=\beta +1$,  we let
$X^{(\alpha)}=(X^{(\beta)})'$ and for limit ordinal $\alpha$ we
let
 $X^{(\alpha)}=\bigcap_{\beta<\alpha}X^{(\beta)}$. The
{\it Cantor-Bendixson rank} of $X$ is the first ordinal
$\alpha<\omega_1$  such that $X^{(\alpha)}=\emptyset$.  The {\it
Cantor-Bendixson rank}, denoted by $r_{cb}(x)$, of $x\in X$ is the
first ordinal $\alpha<\omega_1$  such that $x\in X^{(\alpha)}$ and
$x\notin X^{\alpha+1}$.

\medskip

First, we need to show several auxiliary lemmas.

\medskip

We say that a sequence $(A_n)_{n \in \N}$ of subsets of a space
$X$ converges to a subset $A \subseteq X$, in symbols
$A_n\longrightarrow A$,  if for every $V\in \mathcal{N}(A)$ there
is $m \in \N$ such that $A_n\subseteq V$ for all $m \leq  n$,
$n\in \N$.

\begin{lema} \label{rango 5} Let $(X,f)$ be a dynamical system where $X$ is a compact metric countable space and $x \in X$. Assume that
 that  $r_{cb}(x)=1$  and $f(x)=x$.  If $x_n \to x$ and for every $V \in \mathcal{N}(x)$ and for every $n \in \N$ there is $m \in \N$ such that $O_f(x_k) \cap \big( V \setminus \{ f^i(x_j) : i, j \leq n \}) = \emptyset $ for each $m \leq k \in \N$, then $O_f(x_n) \to x$
\end{lema}

\proof Fix a clopen neighborhood   $V\in \mathcal{N}(x)$.  Since
$r_{cb}(x)=1$, we can assume, without loss of generality, that $V
\setminus\{x\}$ is discrete and that $x_n \in V \setminus \{x\}$
for all $n \in N$. Suppose, towards a contradiction, that
$B=\{n\in\mathbb N: \mathcal O_f(x_n)\not\subseteq V\}$ is
infinite.  Let $n_0  = \min{B}$ and choose $z_{0} \in \mathcal
O_f(x_{n_0}) \cap V_0$ so that $f(z_{0})\notin V$. Suppose that we
have defined $\{ n_i : i < k\} \subseteq B$ and, for each $i < k$,
$z_{i} \in \mathcal O_f(x_{n_i}) \cap V$ such that $z_{i} \neq
z_{j}$, for all $j < i$, and $f(z_{i})\notin V$. For each $i < k$
choose $l_i \in \N$ so that $f^{l_i}(x_{n_i}) = z_i$. Pick a
positive integer  $l > \max\{l_0,\cdots, l_{k-1}\} +
\max\{n_0,\cdots, n_{k-1}\} + 1$. Now, choose $n_k \in B$ so that
$O_f(x_{n_k}) \cap \Big(V \setminus \{ f^i(x_j) : i, j \leq l
\}\Big) = \emptyset$. Then, there is $z_k \in  \mathcal
O_f(x_{n_k}) \cap \big(V \setminus \{ z_i : i < k\}\big)$ so that
$f(z_{k})\notin V$. By construction the set  $\{ z_{n_k} : k \in
\N \}$ is infinite and it is contained in $V$ which implies that
$z_{n_k} \xrightarrow[k \to \infty]{} x$. But this implies that
$f(z_{n_k})\longrightarrow f(x)=x$ which is a contradiction. \finp

\begin{corolario} \label{rango 1}
Let $(X,f)$ be a dynamical system where $X$ is a compact metric
countable space. Assume that $x\in X$ satisfies that $r_{cb}(x)=1$
and $f(x)=x$.  If $x_n \to x$ and $x_n$ is periodic for each $n
\in \N$, then $O_f(x_n) \to x$.
\end{corolario}

The previous corollary is not true if we only assume that each
$x_n$ is eventually periodic. In fact, consider the space $X = \{
\frac{1}{n} : n \in \mathbb{N} \} \cup \{0\}$ and the function
\[
f(x)=\left\{\begin{array}{ll}
                         0, & \hbox{$x=0$} \\
                         1, & \hbox{$x=1$} \\
                        \frac{1}{n} , & \hbox{$x=\frac{1}{n+1}$.}
                       \end{array}
                     \right.
\]
Notice that $1\in \mathcal{O}_f(\frac{1}{n})$, for all $n$. Hence
$\mathcal{O}_f(\frac{1}{n})\not\rightarrow 0$.

\begin{lema} \label{orbitas atrapadas} Let $(X,f)$ be a dynamical system where $X$ is a compact metric space.
 Let $x\in X$   and assume that  $(x_n)_{n\in\mathbb N}$ is a sequence of
points in $X$ such that $x_n\longrightarrow x$. If  $g=f^l$ for
some positive $l \in \N$ and $\mathcal O_g(x_n) \to x$, then
$O_f(x_n) \to O_f(x)$.
\end{lema}

\proof Fix $V\in\mathcal N(\mathcal O_f(x))$ and consider the
following open neighborhood of $x$:
\[
U=\bigcap_{i=0}^{l-1} f^{-i}(V).
\]
Choose $m \in\mathbb N$ so that $O_g(x_n)\subseteq U$, for each
$n\geq m$. By Lemma \ref{descomposici�n de �rbita},   we know
that
$$
\mathcal{O}_f(x_n)=\mathcal{O}_g(x_n)\cup
f[\mathcal{O}_g(x_n)]\cup\ldots\cup f^{l-1}[\mathcal{O}_g(x_n)],
$$
for each $n\in\mathbb N$. Notice that $f^i [U] \subseteq V$ for
all $i\leq l-1$. Hence, if $i\leq l-1$, then
$f^i[\mathcal{O}_g(x_n)]\subseteq V$ for every   $n\geq m$.
Therefore, $\mathcal O_f(x_n)\subseteq V$ for each $n\geq
 m$.
\finp

\begin{lema}\label{orbitas atrapadas 2}
Let $(X,f)$ be a dynamical system such that $X$ is a compact
metric countable space and  every point of $X'$ is  periodic. If
$x\in X$ is a fixed point of $f$ and $(x_n)_{n\in\mathbb N}$ is a
sequence of  periodic points converging to $x$, then $\mathcal
O_f(x_n) \to x$.
\end{lema}

\proof We shall prove this theorem by induction on the
Cantor-Bendixon rank of $x$.  The case when $x$ isolated is
trivial, since eventually we have that $x=x_n$ and $x$ is fixed.
The case $r_{cb}(a)=1$ was already proved in  Corollary \ref{rango
1}. Assume that $r_{cb}(a)=\alpha$ and that the result holds for
all points of $X$ of rank smaller than $\alpha$. Fix a clopen $V
\in\mathcal N(x)$ and assume, without loss of generality, that
$r_{cb}(y) <\alpha$ for each $y \in V\setminus \{x\}$. Also assume
that  $x_n \in V \setminus \{x\}$ for all $n \in \N$. Suppose that
$B=\{n\in\mathbb N: \mathcal O_f(x_n)\not\subseteq V\}$ is
infinite. For each $n\in B$ choose $y_n\in \mathcal O_f(x_n) \cap
V$ so that $f(y_n)\notin V$. As the set $\{ y_n : n \in B \}$ is
infinite, we can find $y \in V$ and an increasing sequence
$(n_k)_{k\in\mathbb N}$ in $B$ such that $y_{n_k}\longrightarrow
y$. It is evident that $f(y) \notin V$ and hence $f(y) \neq y \neq
x$. We claim that $y$ is periodic.   In fact, if $y \in X'$, then
$y$ is periodic by assumption. Thus, suppose that  $y$ is an
isolated point. Then eventually $y_{n_k}= y$ and hence $y \in
\mathcal O_f(y_{n_k})$, which implies that $y$ is periodic. Let $l
=|\mathcal O_f(d)|$ and set $g=f^l$. It is clear that  each
$f$-periodic point is $g$-periodic. Since $y_{n_k} \in
O_f(x_{n_k})$ for every $k \in \N$, then we must have that  each
$y_{n_k}$ is $g$-periodic. As $x$ is a fixed point of $f$, then $
x \notin\mathcal{O}_f(y)$. Hence, there is a clopen  $U
\in\mathcal N(\mathcal O_f(y))$  such that $x \notin U$. We know
that  $r_{cb}(y) <\alpha$. By applying the inductive hypothesis to
$(X,g)$, we obtain that  $\mathcal O_g(y_{n_k})\xrightarrow[k \to
\infty]{} y$. It then follows from Lemma \ref{orbitas atrapadas}
that there exists $i \in \mathbb{N}$ such that $\mathcal
O_f(y_{n_j}) \subseteq U$, for each $i \leq j \in \N$. But this is
impossible since $x_{n_k} \in \mathcal O_f(x_{n_k}) =
\mathcal{O}_f(y_{n_k})$ for each $k \in \N$ and
$x_{n_k}\xrightarrow[k \to \infty]{} x \notin U$. \finp

\begin{corolario}\label{orbitasatrapadas3}
Let $(X,f)$ be a dynamical system such that $X$ is a compact
metric countable space and  every point of $X'$ is  periodic. If
$(x_n)_{n\in\mathbb N}$ is a sequence of periodic points of $X$
converging to $x$, then $\mathcal O_f(x_n) \to \mathcal O_f(x)$.
\end{corolario}

\proof If $x$ is isolated, then  eventually  $x_n=x$ and the
result follows. Suppose that $x\in X'$ and let $l$ be the period
of $x$. Let $g=f^l$, then $g(x)=x$ and by Lemma \ref{orbitas
atrapadas 2} applied to $(X,g)$, we obtain that
$\mathcal{O}_g(x_n)\rightarrow x$. Now, by Lemma \ref{orbitas
atrapadas}, we get the conclusion. \finp

For points with an infinite orbit, we have the following result
which follows directly from Lemma \ref{rango 5}.

\begin{lema}\label{caso rango CB= 1}
Let $(X,f)$ be a dynamical system such that $X$ is a compact
metric countable space. If $x\in X$ has infinite orbit and there
is $y \in \omega_f(x)$  such that $r_{cb}(y) =1$ and $f(y)=y$,
then $f^n(x)\longrightarrow y$.
\end{lema}

\begin{lema}\label{orbita converge para X' periodicos}Let $(X,f)$ be a dynamical system such that $X$ is a compact metric countable space and  every point of $X'$
is  periodic. If $x$ has an infinite orbit and  $y \in w_f(x)$  is
fixed, then $f^n(x)\longrightarrow y$.
\end{lema}

\proof We will prove it by induction on the Cantor-Bendixon rank
of $y$ and for each continuous function $f: X \to X$. The case in
which $r_{cb}(y)=1$ was already proved on lemma \ref{caso rango
CB= 1}. Suppose  $r_{cb}(y)=\alpha > 1$ and that the result holds
for every continuous $g:X\rightarrow X$ and for every point of
$X'$ with
 Cantor-Bendixon rank $<\alpha$.  Choose a clopen $V \in \mathcal N(y)$  such that every point of $V\setminus \{y\}$ has
 Cantor-Bendixon rank $<\alpha$.
  Assume that the set $A = \{n\in \mathbb N: f^n(x)\in
V \;\mbox{and}\;f^{n+1}(x)\notin V\}$ is infinite. Then there is
$z \in V \setminus\{y\}$ and an increasing sequence
$(n_k)_{k\in\mathbb N}$ in $A$ such that $f^{n_k}(x)\xrightarrow[k
\to \infty]{}  z$ and, $f^{n_k}(x) \in V$ and $f^{n_k+1}(x)\notin
V$ for all $k \in \N$. Clearly, $f(z) \neq z \neq y$. Set $l
=|\mathcal O_f(z)|$. Now for each $k\in \mathbb N$ pick $t_k$,
$r_k\in \mathbb N$ so that $n_k=t_kl+r_k$ and $r_k < l$. Choose
$r<l$ such that $B = \{n_k: r_k=r\}$ is infinite and set $g=f^l$.
As $f^{n_k}(x)=g^{t_k}(f^r(x))$ for each $n_k \in B$, then
$z\in\omega_g(w)$ where $w=f^r(x)$. Since  $g(z)=z$ and $r_{cb}(z)
<\alpha$,  by the inductive hypothesis applied to $(X,g)$, it
follows that $g^n(w)\longrightarrow z$. Hence, we have that
$f^{ln}(x)\xrightarrow[n \to \infty]{} z$. According to Lemma
\ref{orbitaomega}, we then have that $\omega_f(x) = \mathcal
O_f(z)$, which is a contradiction since $y \notin
\mathcal{O}_f(z)$ and $y \in \omega_f(x)$. \finp

\begin{teorema}\label{w_f(y)=orbita periodica} Let $(X,f)$ be a dynamical system such that $X$ is a compact metric countable space and  every point of $X'$
is  periodic. For every  $x\in X$, there is a periodic point
$y\in X$ such that $\omega_f(x)=\mathcal O_f(y)$.
\end{teorema}

\proof The case when $\mathcal O_f(x)$ is finite it is evident.
Assume that $\mathcal O_f(x)$ is infinite and let $y \in
\omega_f(x)$. By assumption, $y$ is periodic. Set $l=|\mathcal
O_f(y)|$ and $g=f^l$. Then, we have that $g(y)=y$. Choose an
increasing sequence $(n_k)_{k\in \N}$ for which $f^{n_k}(x)$
converges to $y$. By passing to a subsequence, we assume, without
loss of generality, there is  $i<l$ such that  $n_k=lt_k+i$  where
$t_k \in \N$, for all $k$. Set $z =f^i(x)$. Then, we obtain that
$g^{t_k}(z)\xrightarrow[k \to \infty]{} y$. This implies that $y
\in \mathcal{O}_g(z)$. By Lemma \ref{orbita converge para X'
periodicos}, $g^n(z)\to y$ and, by Lemma \ref{orbitaomega}, obtain
that $\omega_f(x) = \mathcal O_f(y)$. \finp

The two following theorems will allow us to conclude that given a
dynamical system $(X,f)$, where each element of $X'$ is a periodic
point and given  $x\in X$, either $f^q$ is discontinuous at $x$
for all $q\in\mathbb N^\ast$ or  $f^q$ is continuous at $x$ for
all $q\in\mathbb N^\ast$.

\begin{teorema}\label{general-punto fijo}
Let $(X,f)$ be a dynamical system such that $X$ is a compact
metric countable space and  every point of $X'$ is  periodic. If
$x\in X'$ is a fixed point, then either  $f^p$ is continuous at
$x$, for every $p \in\mathbb N^\ast$, or $f^p$ is discontinuous at
$x$, for every $p\in\mathbb N^\ast$.
\end{teorema}

\proof Let $x \in X'$ be a fixed point. Suppose that there exist
$p, q\in\mathbb N^\ast$ such that $f^p$ is continuous at $x$  and
$f^q$ is discontinuous at $x$. By compactness, there is a sequence
$(x_n)_{n\in\mathbb N}$ in $X$  such that $x_n\longrightarrow x$
and $y_n =f^q(x_n)\longrightarrow y \neq x$.   According to
Theorem \ref{w_f(y)=orbita periodica}, for each $n \in \N$ there
is a periodic point $z_n \in X$ such that $\omega_f(x_n)= \mathcal
O_f(z_n)$.  Without loss of generality, we may assume that $z_n
\longrightarrow z$. Clearly $z$ is periodic.  By Corollary
\ref{orbitasatrapadas3}, we obtain that $\mathcal O_f(z_n) \to
\mathcal O_f(z)$. Since $f^q(x_n) \in \omega_f(x_n)= \mathcal
O_f(z_n)$, hence $y\in \mathcal{O}_f(z)$. On the other hand, by
the continuity of $f^p$ at $x$, we must have that
$f^p(x_n)\longrightarrow x$. Since $f^p(x_n) \in \omega_f(x_n)$,
we conclude, as before, that $x\in \mathcal{O}_f(z)$. Since $x$ is
fixed and $y$ is periodic, then $x=y$, which is a contradiction.
\finp

\begin{lema}
\label{orbita periodica discontinua} Let $(X,f)$ be a dynamical
system such that $X$ is a compact metric space and  every point of
$X'$ is  periodic. If  $p\in \mathbb{N}^\ast$ and  $f^p$ is
continuous at the point $x \in X'$, then $f^p$ is continuous at
every point of $\mathcal O_f(x)$.
\end{lema}

\proof Let  $p\in \mathbb{N}^\ast$ and suppose that $f^p$ is
continuous at the point $x \in X'$. Let $n$ be the period of $x$.
Consider the point  $y=f^l(x)$ for some $l<n$. Let $(y_k)_{k \in
\mathbb{N}}$ be a sequence in $X$ such that $y_k\longrightarrow
y$. Then, we have that $f^{(n-l)}(y_k)\longrightarrow x$ and so
$f^p(f^{(n-l)}(y_k))\longrightarrow f^p(x)$.
 Since $p+m=m+p$  for every $m\in \mathbb N$, we must have that
$f^p \circ f^{(n-l)}=f^{(n-l)} \circ f^p$. Thus, we obtain that
$f^n(f^p(y_k))\longrightarrow f^l(f^p(x)) =f^p(f^l(x)) =f^p(y)$.
Now, assume  that $f^p(y_k)\not \longrightarrow f^p(y)$. By
passing to a subsequence, without loss of generality, we may
assume that $f^p(y_k)\longrightarrow z \neq f^p(y)$. Hence,
$f^n(f^p(y_k)) \longrightarrow f^n(z)$ and so $f^n(z) = f^p(y)$.
Since  $f^p(y)\in \mathcal O_f(x)$ and $z$ is periodic, then $z$
has period $n$ and hence $z=f^n(z)=f^p(y)$,  but this is
impossible. Therefore, $f^p$ is continuous at $y$. \finp

Now we are ready to state the main result of this section.

\begin{teorema}
\label{general-punto periodico} Let $(X,f)$ be a dynamical system
such that $X$ is a compact metric countable space and  every point
of $X'$ is  periodic. Then, for each $x \in X$ either $f^p$ is
discontinuous at $x$, for all $p\in\mathbb N^\ast$, or  $f^p$ is
continuous at $x$, for all $p\in\mathbb N^\ast$.
\end{teorema}

\proof Fix $x \in X'$ and $p, q\in\mathbb N^\ast$. Without loss of
generality, assume that $f^p$ is continuous at $x$. According to
Lemma \ref{orbita periodica discontinua}, we know that $f^p$ is
continuous at every point of $\mathcal O_f(x)$.  Now, set
$n=|\mathcal O_f(x)|=$ and $g=f^n$.  By Lemma \ref{iteradas de
potencias}, $g^p = f^p \circ f^n$,  hence we must have that $g^p$
is also continuous at $x \in X$. Notice that each point of
$\mathcal O_f(x)$ is a fixed point of $g$. Then, by Theorem
\ref{general-punto fijo},  $g^q$ is continuous at each point of
$\mathcal O_f(x)$. Choose $l < n$ so that $ q \in(n\mathbb
N+l)^\ast$. By Lemma \ref{iteradas de potencias}, $f^q = g^{q}
\circ f^l$.  So, $f^q$ is also continuous at $x$. \finp

We have already mentioned that for $X$ a convergent sequence, it
was shown  in \cite{GarciaSanchis} that  $f^p$ is continuous for
every $p\in\mathbb{N}^*$ or $f^p$ is discontinuous for every
$p\in\mathbb{N}^*$. In the next theorem we shall extend this
result to any compact metric countable space with finitely many
accumulation points (in the next section we will see that this
result cannot be extended to a space with CB-rank equal to 2).

\medskip

The following lemma will help us to consider only orbits without
isolated points.

\begin{lema}
\label{punto aislado en orbita} If $\mathcal{O}_f(x)$ has an
isolated point, then every $f^{p}$ is continuous at $x$, for every
$p\in \N^*$
\end{lema}

\proof Suppose $f^n(x)$ is isolated. Let $x_k\rightarrow x$. Then
there is $k_0$ such that $f^n(x_k)=f^n(x)$ for all $k\geq k_0$.
Therefore $f^p(x_k)=f^p(x)$ for all  $k\geq k_0$.
\endproof

\medskip

To proof our theorem we need a synchronization  of the
$q$-iterations of points near to an orbit of a periodic point.

\begin{lema}\label{sincronia}
Suppose $X'$ is finite, $z$ periodic with period $s$ and
$\mathcal{O}_f(z)\subseteq X'$. Let $V_j$ be pairwise disjoint
clopen sets, for $0\leq j<s$, such that $V_j\cap X'=\{f^j(z)\}$
for all $j$. Let $p\in \N^*$ be such that $f^p(z)=z$    and $x\in
V_0$. If $f^p(x)\in V_0$, then
\[
\mathcal{O}_f(x)\subseteq^* V_0\cup\cdots\cup V_{s-1}
\]
and $f^q(x)\in V_j$ whenever that $q\in (s\N+j)^*$ with $0\leq
j<s$.
\end{lema}

\proof Since $f^p(x)\in V_0$, then $\{m\in\N:\; f^m(x)\in V_0\}\in
p$. As $f^p(z)=z$, then $s\N\in p$. Therefore
\[
A_0=\{ns:\; f^{ns}(x)\in V_0\}\in p.
\]
Fix  $j < s$. Since $f^j(z)\in V_j$ and $V_j$ has only one
accumulation point, then $f^j[V_0]\subseteq^* V_j$ by the
continuity of $f$. Therefore
\[
A_j=\{ns+j:\; f^{ns+j}(x)\in V_j\}
\]
is infinite for each $j<s$. On the other hand, as $z$ has period
$s$, then  $f^s[V_j]\subseteq^* V_j$. Thus $A_j+s\subseteq^* A_j$
for all $j$. Thus $A_j=^* s\N+j$. In particular, $f^{ns+j}(x)\in
V_j$ for  almost all $n$. This says that
$\mathcal{O}_f(x)\subseteq^* V_0\cup\cdots\cup V_{s-1}$. Moreover
$f^q(x)\in V_j$ provided that  $q\in (s\N+j)^*$ for some $0\leq
j<s$.
\endproof

We are ready to prove the last main theorem of the section.

\begin{teorema}
\label{Rango finito} Suppose $X'$ is finite. If $f^p$ is
continuous,  for some $p\in\N^*$, then $f^q$ is continuous for all
$q\in\N^*$.
\end{teorema}

The proof will follow from the next  lemmas.

\medskip

We omit the proof of the following easy lemma which takes care of
the case when there is an isolated points in an orbit.

\begin{lema}\label{aislado} Let $z_k \in X$, for each $k \in \N$, and let $z, u \in X$ be two  periodic points such that $\mathcal{O}_f(z_k) \subseteq^* \mathcal{O}_f(u)$ for every $k \in \N$. If $f^p(z_k) \to z$, then $\mathcal{O}_f(u) = \mathcal{O}_f(z)$.
\end{lema}

\begin{lema}
\label{continuidad-periodicos} Suppose $X'$ is finite, $z\in X$
periodic and $\mathcal{O}_f(z)\subseteq X'$. If $f^p$ is
continuous at $z$,  for some $p\in\N^*$, then $f^q$ is continuous
at $z$,  for all $q\in\N^*$.
\end{lema}

\proof  Let $s$ be the period of $z$ and let  $V_j$, for $0\leq
j<s$, be clopen sets as in the hypothesis of lemma
\ref{sincronia}. Notice that if $f^p(x)=f^i(z)$ with $0\leq i<s$
and $f^p$ is continuous at $z$, then taking $r=p+s-i$, we have
that $f^r(z)=z$ and $f^r$ is continuous at $z$. Therefore, without
loss of generality, we will assume that $f^p(z)=z$. Fix $q\in\N^*$
and choose $0\leq j<s$ so that $q\in (s\N+j)^*$. Let $(z_k)_{k \in
\N}$ be a sequence in $V_0$ converging to $ z$.  Since $f^p(z_k)
\to z \in V_0$, by lemma \ref{sincronia}, $f^q(z_k)\in V_j$ for
all  $k$. By Lemma \ref{piteradaperiodico2} , we know that
$f^j(z)=f^q(z)$.

We claim that $f^q(z_k)$ converges to $f^j(z)$. Otherwise,  there
is an isolate point $u$ such that  $f^q(z_k)=u$ for infinitely
many $k$. Hence, $u$ is periodic and $\mathcal{O}_f(z_k)
\subseteq^* \mathcal{O}_f(u)$ for infinitely many $k \in \N$. By
Lemma \ref{aislado}, $u \in \mathcal{O}_f(z)\subseteq X'$ which is
a contradiction. Therefore, $f^q$ is continuous at $z$.
\endproof

\begin{lema}
Suppose $X'$ is finite, $z\in X$ periodic and
$\mathcal{O}_f(z)\subseteq X'$. Let $x\in X'\setminus
\mathcal{O}_f(z)$ be  such that $f^i(x) \in \mathcal{O}_f(z) $ for
some $i \in \N$. If $f^p$ is continuous, for some $p\in\N^*$, then
$f^q$ is continuous at $x$,  for all $q\in\N^*$.
\end{lema}

\proof Let $s$ be the period of $z$ and   $V_j$, for $0\leq j<s$,
be clopen sets as in the hypothesis of lemma \ref{sincronia} and
put  $V=  V_0\cup\cdots\cup V_{s-1}$. Fix  $q\in\N^*$ and choose
$0\leq j<s$ so that  $q\in (s\N+j)^*$. Suppose that
$x_k\rightarrow x$.  We can assume, without loss of generality,
that  $i$ is the smallest $n$ such that $f^n(x)$  is periodic,
$f^i(x)=z$ and $f^i(x_k)\in V_0$ for all $k$. Suppose
$f^q(x_k)\rightarrow w$. We will show that $w=f^q(x)$. We claim
that $w$ is a limit point. Suppose that $w$ is isolated. As in the
proof of the previous lemma, we obtain that  $w$ is periodic and
$\mathcal{O}_f(x_k) \subseteq \mathcal{O}_f(w)$ for infinitely
many $k \in \N$. As $f^p(x_k) \to f^p(x) \in \mathcal{O}_f(z)$, by
Lemma \ref{aislado}, $w \in \mathcal{O}_f(w) = \mathcal{O}_f(z)$
which is impossible. Therefore, $w$ must be a limit point. Notice
that by lemma \ref{sincronia}, $\mathcal{O}_f(f^i(x_k))\subseteq^*
V$  for large enough $k$. Thus $w\in V$ and, being non isolated,
it belongs to the orbit of $z$. On the other hand, since $f^p$ is
continuous, then by lemma \ref{continuidad-periodicos}, $f^q$ is
continuous at $f^i(x)$ and thus  $f^q(f^i(x_k))\rightarrow
f^q(f^i(x))$. Since $f^i(f^q(x_k)) \rightarrow f^i(w)$ and
$f^q\circ f^i=f^i\circ f^q$, then $f^i(w)=f^i(f^q(x))$.   As
$f^q(x)$ and $w$ are both  in the orbit of $z$, then necessarily
$w=f^q(x)$.
\endproof

%%%%%%%%%%%%%%%%%%%%%%%%%%%%%%%%%%%%%%%%%%%%%%%%%%%%%%%%%%%%%%%%%%%%%%%%%%%%%%%%%%%%%%%%%%%%%%%%%%%%%%%%%%%%%%%%%%%%%%%
\section{An example}
%%%%%%%%%%%%%%%%%%%%%%

We construct a  dynamical system $(X,f)$ where $X$ is a compact
metric countable space,  the orbit of each accumulation point is
finite and that there are $p$, $q\in\mathbb N^\ast$ such that
$f^p$ is continuous on $X$ and $f^q$ is discontinuous at some
point of $X$. This shows that the hypothesis ``every accumulation
point is periodic'' of Theorem \ref{general-punto periodico}
cannot be weakened to just asking that the orbit of such points
are finite.

\bigskip

\begin{example}\label{primero} Consider the countable ordinal space  $\omega^2 +1$ which will be identified with a suitable
subspace of $\mathbb{R}$:
$$
X= \Big(\bigcup_{m\in\mathbb N}D_m \Big) \cup \{d_m : m \in \N
\}\cup\{d\},
$$
where $(d_m)_{m \in \N}$ is a strictly increasing  sequence
converging to $d$, $D_m=\{d^m_n:n\in\mathbb N\}$ is a strictly
increasing sequence contained in $(d_{m-1},d_{m})$ converging to
$d_m$, for each $m \in \N \setminus \{1\}$, and
$D_1=\{d^1_n:n\in\mathbb N\}$ is a strictly increasing sequence
contained in $(-\infty,d_{1})$ converging to $d_1$. For notational
convenience, we shall assume that $0 \notin \mathbb{N}$ and
$\mathbb P$ stands for the set of prime numbers. Now, we define
the function $f:X\rightarrow X$ as follows:

\begin{enumerate}
\item[(i)] $f(d)=d$, $f(d_1)=d$ and $f(d_m)=d_{m-1}$ for each
$m>1$.

\item[(ii)] $f(d^1_n)=d_n^n$, for each $n\in\mathbb P$.

\item[(iii)]  $f(d^1_n)=d$ for each $n\not\in\mathbb P$.

\item[(iv)]  $f(d^m_n)=d^{m-1}_n$ whenever $m> 1$ and $n\not\in
   \mathbb P$.

\item[(vi)]  $f(d^m_n)=d^{m-1}_n$ whenever $1<m\leq n$ and
$n\in\mathbb P$.

\item[(vii)] $f(d^m_n)= d^{m-1}_{n-1}$ whenever $m>n$ and
$n\in\mathbb P$.
\end{enumerate}
It is not hard to prove that $f$ is continuous and notice that
every point is eventually periodic. The required dynamical system
will be $(X,f)$. We have the following consequences directly from
the definition:
\begin{enumerate}
\item $d_m$ is eventually periodic and $\mathcal
O_f(d_m)=\{d_m,d_{m-1},d_{m-2},\cdots,d_{1},d\}$ for each $m \in
\N$.

\item $d^1_n$ is eventually periodic and $\mathcal
O_f(d^1_n)=\{d^1_n,d\}$ for  each $n\not\in \mathbb P$.

\item $d^n_n$ is periodic and $\mathcal
O_f(d_n^n)=\{d_n^n,d_n^{n-1},d_n^{n-2},\cdots,d_n^{1}\}$ for each
$n\in\mathbb P$.

\item $d^m_n$ is eventually periodic and $\mathcal
O_f(d^m_n)=\{d^m_n,d^{m-1}_{n},d_{n}^{m-2},\cdots,d_{n}^{1},d\}$
provided that  $n\not\in \mathbb P$ and $m>1$.

\item $d^m_n$ is eventually periodic and $\mathcal
O_f(d^m_n)=\{d^m_n,d^{m-1}_{n-1},d_{n-1}^{m-2},\cdots,d_{n-1}^{1},d\}$
provided that  $n < m$ and $n\in\mathbb P$.

\item $d^m_n$ is periodic and $\mathcal
O_f(d^m_n)=\{d^m_n,d^{m-1}_{n},d_{n}^{m-2},\cdots,d_{n}^{1},d_n^n,d_n^{n-1},d_n^{n-2},\cdots,d_n^{m+1}\}$
provided that  $1 < m < n$ and $n\in\mathbb P$.
\end{enumerate}
Hence, we obtain that:
\begin{enumerate}
\item[(a)]  $f[D_m]\subseteq D_{m-1}$ for all $m>1$.

\item[(b)]  $f[D_1\setminus\{d^1_n:n\in\mathbb P\}]=\{d\}$ and
$f[\{d^1_n:n\in\mathbb P\}]=\{d_n^n:n\in\mathbb P\}$.

\item[(c)] Fore each $x\notin \bigcup_{n\in\mathbb P}\mathcal
O_f(d_n^n)$ there exists $l \in\mathbb N$ such that $f^l(x)=d$.
\end{enumerate}
\end{example}

To analyze the behavior of the $p$-iterates of $f$, we shall need
some preliminary lemmas.

\medskip

From the computation of the orbits given above and lemma
\ref{piteradaperiodico2} we have the following result.

\begin{lema}\label{ejemplo-p-iteradas} Let $p\in\mathbb N^*$ and let $l_n \in \N$ (depending on
$p$) be such that $0\leq l_n <n$ and $p \in (n\mathbb{N}+l_n)^*$.
Then, we have that
\begin{itemize}
\item[(i)] $f^p(d^m_n)=d$ whenever  $n\in\mathbb P$ and $m>n$.

\item[(ii)] $f^p(d^m_n)=d$ whenever  $n\not\in \mathbb P$ and
$m>1$.

\item[(iii)] $f^p(d^m_n)=f^{l_n}(d^m_n)$, when  $n\in \mathbb{P}$
and $m\leq n$. In particular, $f^p(d^m_n)=d^{n-(l_n-m)}_n$
provided that $m\leq l_n< n \in \mathbb{P}$.
\end{itemize}
\end{lema}

\begin{teorema}\label{mainresult} Let $(X,f)$ be the dynamical system constructed
above.
\begin{enumerate}
\item If  $p\in \bigcap_{n\in\mathbb P}(n\mathbb N+(n-1))^\ast$,
then $f^p$ is discontinuous at $d$.

\item If $p\in \bigcap_{n\in\mathbb P}(n\mathbb N+
\frac{n+1}{2})^\ast$, then $f^p$ is continuous on $X$.
\end{enumerate}
\end{teorema}

\proof $(1).$ Let $p\in \bigcap_{n\in\mathbb P}(n\mathbb
N+(n-1))^\ast$. According to  Lemma \ref{ejemplo-p-iteradas}(iii),
we know that $f^p(d_n^n)= f^{n-1}(d_n^n)= d^1_n$ for all
$n\in\mathbb P$. Hence, we obtain that  the sequence
$(f^p(d_n^n))_{n\in \mathbb P}$ converges to $d_1$, but the
sequence $(d_n^n)_{n\in \mathbb P}$ converges to $d$ and $f(d)=d$.
Therefore, $f^p$ is not continuous at $d$.

$(2).$ Let $p\in\bigcap_{n\in\mathbb P}(n\mathbb
N+\frac{n+1}{2})^\ast$.  We first show that $f^p$ is continuous at
$d_m$ for every $m \in \N$. Let $m \in \N$ and assume $x_k
\xrightarrow[k \to \infty]{} d_m$. We remark that $f^p(d_m)=d$.
Without loss of generality, suppose that $x_k=d^m_{n_k}$ where  $m
< \frac{n_k+1}{2}$ for each $k \in \N$.  If  $n_k \not\in \mathbb
P$ for some $k \in \N$,  by Lemma \ref{ejemplo-p-iteradas}(ii),
then we have that $f^p(x_k)=d$. Thus, we may suppose that $n_k \in
\mathbb{P}$ for all $k \in \N$. It then follows from Lemma
\ref{ejemplo-p-iteradas}(iii) that
$$
f^p(x_k) = f^p(d^m_{n_k}) =  d_{n_k}^{m + \frac{n_k-1}{2}}
\xrightarrow[k \to \infty]{} d = f^p(d_m).
$$
Next, we shall show that $f^p$ is continuous at $d$.  To do that
let us assume that  $x_k \longrightarrow d$ and $\{ x_k : k \in \N
\} \cap \{ d_m : m \in \N\} = \emptyset$. Write
$x_k=d^{m_k}_{n_k}$ for each $k \in \N$ . As above, without loss
of generality, we may suppose that $n_k\in \mathbb{P}$ for each $k
\in \N$. If $n_k <m_k$ for some $k \in \N$, by Lemma
\ref{ejemplo-p-iteradas}(i), then $f^p(d^{m_k}_{n_k}) =  d$.
Thus, we can also assume that  $m_k\leq n_k$ for all $k \in \N$.
In virtue of lemma \ref{ejemplo-p-iteradas}(iii), we have that
$$
f^p(x_k) = f^p(d^{m_k}_{n_k}) =  d_{n_k}^{n_k - \frac{n_k+1}{2} +
m_k} = d_{n_k}^{m_k + \frac{n_k-1}{2}} \xrightarrow[k \to
\infty]{} d = f^p(d_m).
$$
\finp

To finish our task, we  shall show that there are ultrafilters
satisfying  the hypothesis of  Lemma \ref{mainresult}.

\begin{lema}
\label{piteradasprimosrelativos} Let $(n_k)_{k\in \mathbb{N}}$ be
an increasing sequence of pairwise relatively prime natural
numbers. For every sequence $(d_k)_{k \in \N}$ satisfying  $0 \leq
d_k < n_k$ for each $k \in \N$, we have that
$$
\bigcap_{k \in \N} (n_k\mathbb N+d_k)^\ast \neq \emptyset.
$$
\end{lema}

\proof  It is enough to show that the family $\{(n_k\mathbb
N+d_k)^\ast:k\in\mathbb N\}$ has the infinite finite intersection
property. Indeed, let ${k_1},\cdots, {k_l}$ in $\mathbb N$.
Consider the equations system:
$$
\left\{\begin{array}{ll}
                         d_{k_1}\equiv x, & \hbox{$(mod\,\, n_{k_1})$;} \\
                         \vdots & \hbox{} \\
                         d_{k_n}\equiv x, & \hbox{$(mod\,\, n_{k_l})$.}
                       \end{array}
                     \right.
$$
Since the natural numbers $n_{k_1}$, $\ldots$, $n_{k_l}$ are
relatively prime, by the Chinese Remainder Theorem,  this system
has infinitely many solutions. Therefore, the intersection
$$
(n_{k_1}\mathbb N+d_{k_1})\cap\ldots\cap(n_{k_l}\mathbb N+d_{k_n})
$$
is infinite. Therefore,   $\bigcap_{k\in\mathbb N}(n_k\mathbb
N+d_k)^\ast \neq \emptyset$. \finp

\bigskip

We finish this section with  some questions: Given  $p,
q\in\mathbb{N}^*$ such that $p+n \neq q$, for all $n \in
\mathbb{N}$, is there a dynamical system $(X,f)$ and a point $x
\in X$ such that $X$ is a compact metric space, $f^p$ is
continuous at $x$ and $f^q$ is discontinuous at $x$? We remark
that the continuity of $f^p$, for $p \in \N^*$, implies the
continuity of $f^{p+n}$ for each $n \in \mathbb{N}$. In the light
of the example we presented, we would like to know the answer of
previous question for the space $\omega^2 +1$.

\medskip

We thank the referee for providing constructive comments that help
improving the contents of this paper. Mainly, we thank him/her for
pointing out a mistake in the proof of Theorem \ref{Rango finito}.

%%%%%%%%%%%%%%%%%%%%%%%%%%%%%%%%%%%%%%%%%%%%%%%%%%%%%%%%%%%%%%%%%%%%%%%%%
%%%%%%%%%%%%%%%%%%%%%%%%%%%%%%%%%%%%%%%%%%%%%
\bibliographystyle{amsplain}

\end{document}